\documentclass[12pt]{article}
\usepackage{amsmath, amsthm, amssymb} 

\newtheoremstyle{theorem}
  {10pt}		  
  {10pt}  
  {\sl}  
  {\parindent}     
  {\bf}  
  {. }    
  { }    
  {}     
\theoremstyle{theorem}
\newtheorem{theorem}{Theorem}

\newtheorem{remark}[theorem]{Remark}
\newtheorem{example}[theorem]{Example}

\newtheoremstyle{defi}
  {10pt}		  
  {10pt}  
  {\rm}  
  {\parindent}     
  {\bf}  
  {. }    
  { }    
  {}     
\theoremstyle{defi}

\def\proofname{\indent {\sl Proof.}}

\def\ds{{\mathrm{d}}s}
\def\dt{{\mathrm{d}}t}

\usepackage[paperwidth=8.5in,paperheight=11.0in,
  left=1.5in,right=1.0in,top=1.0in,bottom=1.0in]{geometry}

\begin{document}


\author{
{\bf Robert Vrabel}           
\vspace{1mm}\\
{Slovak University of Technology in Bratislava}\\
{Faculty of Materials Science and Technology}\\ 
{J.~Bottu 25, 917 01 Trnava, SLOVAKIA}\\
{robert.vrabel@stuba.sk}\\
}

\title{Singularly perturbed linear Neumann problem with the characteristic roots on the imaginary axis. \\ A non-resonant case\footnote{So far unpublished manuscript was originally written in 2010.}}

\maketitle

\begin{abstract}
In this note we are dealing with the problem of existence and asymptotic behavior of solutions for the non-resonant singularly perturbed linear Neumann boundary value problem 
\begin{eqnarray*}
\epsilon y''+ky=f(t),\quad  k>0,\quad 0<\epsilon<<1,\quad t\in\langle a,b \rangle
\end{eqnarray*}
\begin{equation*} 
y'(a)=0,\quad y'(b)=0.
\end{equation*}
Our approach is based on the analysis of an integral equation equivalent to this problem.

\centerline{}


\centerline{}

\end{abstract}
\section{Introduction}

In this brief paper we will study the singularly perturbed linear problem
\begin{equation}\label{def_DE}
\epsilon y''+ky=f(t),\quad  k>0,\quad 0<\epsilon<<1,\quad f\in C^3\left(\langle a,b\rangle\right)
\end{equation}
with the Neumann boundary condition
\begin{equation}\label{def_BC}
y'(a)=0,\quad y'(b)=0.
\end{equation}
The analysis of the problem under consideration is complicated by the fact that characteristic equation of this differential equation has roots on the imaginary axis, that is, the system is not hyperbolic. For the singularly perturbed dynamical systems the dynamics near a normally-hyperbolic critical manifold is well--known, see e.~g.~\cite{Jo} for geometric approach to the singular perturbation theory and \cite{VraN} for the lower and upper solution method, but if the condition of normal hyperbolicity of a critical manifold is not fulfilled then the problem of existence and asymptotic behavior of solutions is hard solvable in general and leads  to the principal technical difficulties in the nonlinear cases, see for example \cite{Vra}. Thus, the considerations below may be instructive for an analysis of this class of problems. 

Despite these difficulties we will prove, that there exists infinitely many sequences $\left\{\epsilon_n\right\}_{n=0}^\infty,$ $\epsilon_n\rightarrow 0^+$ such that $y_{\epsilon_n}(t)$ converge uniformly to $u(t)$ on $\langle a,b\rangle$ for $\epsilon_n\rightarrow 0^+,$ where $y_{\epsilon_n}$ is a solution of problem (\ref{def_DE}), (\ref{def_BC}) with $\epsilon=\epsilon_n$ and $u$ representing the critical manifold for our system is a solution of the reduced problem $ky=f(t)$ obtained from (\ref{def_DE}) for $\epsilon=0.$ 

Henceforth in the paper we consider for the values of parameter $\epsilon$ only the sets $J_n$ defined as 
\[
J_{n}=\left\langle k\left(\frac{b-a}{(n+1)\pi-\lambda}\right)^2,k\left(\frac{b-a}{n\pi+\lambda}\right)^2\right\rangle,\quad n=0,1,2,\dots,
\]
where  $\lambda>0$ is an arbitrarily small but fixed constant, which guarantees the existence and uniqueness of solutions of (\ref{def_DE}), (\ref{def_BC}), that is, a non-resonant case.
\begin{example}\label{example1}
\rm
Consider the linear problem
\begin{equation*}
\epsilon y''+ky=e^t,\quad t\in\langle a,b \rangle, \quad k>0,\quad 0<\epsilon<<1
\end{equation*}
\begin{equation*}
y'(a)=0,\quad y'(b)=0
\end{equation*}
and its solution
\[
y_\epsilon (t)=\frac{-e^a\cos\left[\sqrt{\frac{k}{\epsilon}}(b-t)\right]+e^b\cos\left[\sqrt{\frac{k}{\epsilon}}(t-a)\right]}{\sqrt{\frac{k}{\epsilon}}(k+\epsilon)\sin\left[\sqrt{\frac{k}{\epsilon}}(b-a)\right]}+\frac{e^t}{k+\epsilon}.
\]
Hence, for every sequence $\left\{\epsilon_n\right\}_{n=0}^\infty,$ $\epsilon_n\in J_n$ the solution of the problem under consideration 
\[
y_{\epsilon_n} (t)=\frac{e^t}{k+\epsilon_n}+{\cal O}(\sqrt\epsilon_n)
\] 
and thus the solutions converge uniformly on the interval $\langle a,b\rangle$ to the solution $u(t)=e^t/k$ of the reduced problem for $n\rightarrow \infty.$
\end{example}

The main result of this note is the following theorem generalizing the Example~\ref{example1} to all right-hand sides $f(t).$  

\section{Main result}

\begin{theorem}\label{main}
For all $f\in C^3\left(\langle a,b\rangle\right)$ and for every sequence $\left\{\epsilon_n\right\}_{n=0}^\infty,$ $\epsilon_n\in J_n$ there exists a unique sequence of the solutions $\left\{y_{\epsilon_n}\right\}_{n=0}^\infty$ of the problem (\ref{def_DE}), (\ref{def_BC}) satisfying 
\begin{equation*}
y_{\epsilon_n}\rightarrow u\ \mathrm{uniformly\ on}\ \langle a,b\rangle\ \mathrm{for}\  n\rightarrow \infty.
\end{equation*}
More precisely,
\begin{equation*}
y_{\epsilon_n}(t)=\frac{f(t)}k+{\cal O}\left(\sqrt{\epsilon_n}\right)\ \mathrm{on}\ \langle a,b\rangle. 
\end{equation*}
\end{theorem}
\proofname\ As first we show that the function
{\setlength\arraycolsep{1pt}
\begin{eqnarray}\label{y}
y_\epsilon(t)&=&\frac{\cos\left[\sqrt{\frac {k}\epsilon}(t-a)\right]\int\limits_{a}^{b}\cos\left[\sqrt{\frac   {k}\epsilon}(b-s)\right]\frac{f(s)}{\epsilon}\ds}{\sqrt{\frac {k}\epsilon}\sin\left[\sqrt{\frac {k}\epsilon}(b-a)\right]}\nonumber\\
&+&\int\limits_{a}^{t}\frac{\sin\left[\sqrt{\frac {k}\epsilon}(t-s)\right]\frac{f(s)}{\epsilon}}{\sqrt{\frac {k}\epsilon}}\ds
\end{eqnarray}}
is a solution of (\ref{def_DE}), (\ref{def_BC}). Differentiating  (\ref{y}) twice, taking into consideration the equality
\begin{equation*}
\frac{\mathrm{d}}{\dt}\int\limits_{a}^{t}H(t,s)f(s)\ds=\int\limits_{a}^{t}\frac{\partial H(t,s)}{\partial t}f(s)\ds+H(t,t)f(t) 
\end{equation*}
we obtain that
{\setlength\arraycolsep{1pt}
\begin{eqnarray}
y'_\epsilon(t)&=&-\frac{\sqrt{\frac {k}\epsilon}\sin\left[\sqrt{\frac {k}\epsilon}(t-a)\right]\int\limits_{a}^{b}\cos\left[\sqrt{\frac{k}\epsilon}(b-s)\right]\frac{f(s)}{\epsilon}\ds}{\sqrt{\frac {k}\epsilon}\sin\left[\sqrt{\frac {k}\epsilon}(b-a)\right]}\nonumber\\
&+&\int\limits_{a}^{t}\frac{\sqrt{\frac {k}\epsilon}\cos\left[\sqrt{\frac {k}\epsilon}(t-s)\right]\frac{f(s)}{\epsilon}}{\sqrt{\frac {k}\epsilon}}\ds,\label{yy}\\
y''_\epsilon(t)&=&-\frac{\left(\sqrt{\frac {k}\epsilon}\right)^2\cos\left[\sqrt{\frac {k}\epsilon}(t-a)\right]\int\limits_{a}^{b}\cos\left[\sqrt{\frac   {k}\epsilon}(b-s)\right]\frac{f(s)}{\epsilon}\ds}{\sqrt{\frac {k}\epsilon}\sin\left[\sqrt{\frac {k}\epsilon}(b-a)\right]}\nonumber\\
&-&\int\limits_{a}^{t}\frac{\left(\sqrt{\frac {k}\epsilon}\right)^2\sin\left[\sqrt{\frac {k}\epsilon}(t-s)\right]\frac{f(s)}{\epsilon}}{\sqrt{\frac {k}\epsilon}}\ds+\frac{f(t)}{\epsilon}.\label{yyy}
\end{eqnarray}}
From (\ref{yyy}) and (\ref{y}), after little algebraic arrangement we get
\[
y''_\epsilon=\frac{k}{\epsilon}\left(-y_\epsilon\right)+\frac{f(t)}{\epsilon},
\]
that is, $y_\epsilon$ is a solution of differential equation (\ref{def_DE}) and from (\ref{yy}) it is easy to verify that this solution satisfies (\ref{def_BC}).

Let $t_0\in\langle a,b\rangle$ is arbitrary but fixed. Let us denote by $I_1$ and $I_2$ the integrals
\begin{equation*}
I_1=:\int\limits_{a}^{b}\cos\left[\sqrt{\frac {k}\epsilon}(b-s)\right]\frac{f(s)}{\epsilon}\ds
\end{equation*}
and
\begin{equation*}
I_2=:\int\limits_{a}^{t_0}\sin\left[\sqrt{\frac {k}\epsilon}\left(t_0-s\right)\right]\frac{f(s)}{\epsilon}\ds.
\end{equation*}
Then
\begin{equation*}
y_\epsilon\left(t_0\right)=\frac{\cos\left[\sqrt{\frac {k}\epsilon}\left(t_0-a\right)\right]I_1}{\sqrt{\frac {k}\epsilon}\sin\left[\sqrt{\frac {k}\epsilon}(b-a)\right]}
+\frac{I_2}{\sqrt{\frac {k}\epsilon}}.
\end{equation*}
 Integrating $I_1$ and $I_2$ by parts we obtain that
{\setlength\arraycolsep{1pt}
\begin{eqnarray*}
I_1&=&\left| \begin{array}{ll}
    h'=\cos\left[\sqrt{\frac {k}\epsilon}(b-s)\right]     & g=\frac{f(s)}{\epsilon}  \\
    h=-\sqrt{\frac {\epsilon}k}\sin\left[\sqrt{\frac {k}\epsilon}(b-s)\right]  & g'=\frac{f'(s)}{\epsilon}
    \end{array} \right|=                  \\
   &=&\sqrt{\frac {\epsilon}k}\sin\left[\sqrt{\frac {k}\epsilon}(b-a)\right]\frac{f(a)}{\epsilon}+ \int\limits_{a}^{b}\sqrt{\frac {\epsilon}k}\sin\left[\sqrt{\frac{k}\epsilon}(b-s)\right]\frac{f'(s)}{\epsilon}\ds,\\
I_2&=&\left| \begin{array}{ll}
    h'=\sin\left[\sqrt{\frac {k}\epsilon}\left(t_0-s\right)\right]     & g=\frac{f(s)}{\epsilon}  \\
    h=\sqrt{\frac {\epsilon}k}\cos\left[\sqrt{\frac {k}\epsilon}\left(t_0-s\right)\right]  & g'=\frac{f'(s)}{\epsilon}
    \end{array} \right|=                  \\
   &=&\frac{\sqrt{\frac {\epsilon}k}f\left(t_0\right)}{\epsilon}-\sqrt{\frac {\epsilon}k}\cos\left[\sqrt{\frac {k}\epsilon}\left(t_0-a\right)\right]\frac{f(a)}{\epsilon}-\int\limits_{a}^{t_0}\sqrt{\frac {\epsilon}k}\cos\left[\sqrt{\frac {k}\epsilon}\left(t_0-s\right)\right]\frac{f'(s)}{\epsilon}\ds.
\end{eqnarray*}}
Thus
{\setlength\arraycolsep{1pt}
\begin{eqnarray*}\label{yto}
y_\epsilon\left(t_0\right)&=&\frac{f\left(t_0\right)}{k}+\frac{\cos\left[\sqrt{\frac {k}\epsilon}\left(t_0-a\right)\right]}{\sin\left[\sqrt{\frac {k}\epsilon}(b-a)\right]}\int\limits_{a}^{b}\sin\left[\sqrt{\frac {k}\epsilon}(b-s)\right]\frac{f'(s)}{k}\ds\nonumber\\
&-&\int\limits_{a}^{t_0}\cos\left[\sqrt{\frac {k}\epsilon}\left(t_0-s\right)\right]\frac{f'(s)}{k}\ds.
\end{eqnarray*}}

Now we estimate the difference $y_\epsilon\left(t_0\right)-\frac{f\left(t_0\right)}{k}.$ We obtain
{\setlength\arraycolsep{1pt}
\begin{eqnarray}\label{est}
\left\vert y_\epsilon\left(t_0\right)-\frac{f\left(t_0\right)}{k}\right\vert&\leq&\frac{1}{k\sin\lambda}\left\vert\int\limits_{a}^{b}\sin\left[\sqrt{\frac {k}\epsilon}(b-s)\right]f'(s)\ds\right\vert\nonumber\\
&+&\frac1k\left\vert\int\limits_{a}^{t_0}\cos\left[\sqrt{\frac {k}\epsilon}\left(t_0-s\right)\right]f'(s)\ds\right\vert.
\end{eqnarray}}
The integrals in  (\ref{est}) converge to zero for $\epsilon=\epsilon_n\in J_n$ as $n\rightarrow \infty.$ Indeed, with respect to the assumption imposed on $f$ we may integrate by parts in (\ref{est}). Thus,
{\setlength\arraycolsep{1pt}
\begin{eqnarray}
&&\int\limits_{a}^{b}\sin\left[\sqrt{\frac {k}\epsilon}(b-s)\right]f'(s)\ds=\left[\sqrt{\frac {\epsilon}k}\cos\left[\sqrt{\frac {k}\epsilon}(b-s)\right]f'(s)\right]^b_a\nonumber\\
&-&\int\limits_{a}^{b}\sqrt{\frac {\epsilon}k}\cos\left[\sqrt{\frac {k}\epsilon}(b-s)\right]f''(s)\ds\nonumber\\
&\leq&\sqrt{\frac {\epsilon}k}\left(\left\vert f'(a)\right\vert+\left\vert f'(b)\right\vert+\left\vert\int\limits_{a}^{b}\cos\left[\sqrt{\frac {k}\epsilon}(b-s)\right]f''(s)\ds\right\vert\ \right)\nonumber\\
&\leq&\sqrt{\frac {\epsilon}k}\left\{\left\vert f'(a)\right\vert+\left\vert f'(b)\right\vert+\sqrt{\frac {\epsilon}k}\left(\left\vert f''(a)\right\vert+\mu_2(b-a)\right)\right\}\label{est2}
\end{eqnarray}}
and
{\setlength\arraycolsep{1pt}
\begin{eqnarray}
&&\int\limits_{a}^{t_0}\cos\left[\sqrt{\frac {k}\epsilon}\left(t_0-s\right)\right]f'(s)\ds=\left[-\sqrt{\frac {\epsilon}k}\sin\left[\sqrt{\frac {k}\epsilon}\left(t_0-s\right)\right]f'(s)\right]^{t_0}_a\nonumber\\
&+&\int\limits_{a}^{t_0}\sqrt{\frac {\epsilon}k}\sin\left[\sqrt{\frac {k}\epsilon}\left(t_0-s\right)\right]f''(s)\ds\nonumber\\
&\leq&\sqrt{\frac {\epsilon}k}\left(\left\vert f'(a)\right\vert+\left\vert\int\limits_{a}^{t_0}\sin\left[\sqrt{\frac {k}\epsilon}\left(t_0-s\right)\right]f''(s)\ds\right\vert\ \right)\nonumber\\
&\leq&\sqrt{\frac {\epsilon}k}\left\{\left\vert f'(a)\right\vert+\sqrt{\frac {\epsilon}k}\left(\mu_1+\left\vert f''(a)\right\vert+\mu_2(b-a)\right)\right\},\label{est3}
\end{eqnarray}}
where $\mu_1=\sup\limits_{t\in\langle a,b\rangle}\left\vert f''(t)\right\vert$ and $\mu_2=\sup\limits_{t\in\langle a,b\rangle}\left\vert f'''(t)\right\vert.$

Substituting (\ref{est2}) and (\ref{est3}) into (\ref{est}) we obtain the a priori estimate of solutions of the problem (\ref{def_DE}), (\ref{def_BC}) for all $t_0\in\langle a,b\rangle$ in the form
{\setlength\arraycolsep{1pt}
\begin{eqnarray}
&&\left\vert y_\epsilon\left(t_0\right)-\frac{f\left(t_0\right)}{k}\right\vert\nonumber\\
&\leq&\frac1{k\sin\lambda}\sqrt{\frac {\epsilon}k}\left\{\left\vert f'(a)\right\vert+\left\vert f'(b)\right\vert+\sqrt{\frac {\epsilon}k}\left(\left\vert f''(a)\right\vert+\mu_2(b-a)\right)\right\}\nonumber\\
&+&\frac1k\sqrt{\frac {\epsilon}k}\left\{\left\vert f'(a)\right\vert+\sqrt{\frac {\epsilon}k}\left(\mu_1+\left\vert f''(a)\right\vert+\mu_2(b-a)\right)\right\}.                   \label{estfinal}
\end{eqnarray}}
 
Because the right-hand side of the inequality (\ref{estfinal}) is independent on $t_0,$ the convergence is uniform on $\langle a,b\rangle.$ Theorem~\ref{main} is proved.

\begin{remark}
\rm
We conclude that in the case that $f'(a)=f'(b)=0,$ that is, the solution $u=f(t)/k$ of a reduced problem satisfies the prescribed boundary conditions (\ref{def_BC}), the convergence rate of the solutions of (\ref{def_DE}), (\ref{def_BC}) to the function $u$ on the interval $\langle a,b \rangle$ is even faster; namely, as ${\cal O}\left(\epsilon_n\right)$ for $\epsilon_n\in J_n$ as follows from (\ref{estfinal}).
\end{remark}

\end{document}